\documentclass{amsart}

\usepackage{amsmath,amsthm,amssymb,amsfonts}

\allowdisplaybreaks

\theoremstyle{plain}
\newtheorem{theo+}           {Theorem}      [section]
\newtheorem{prop+}  [theo+]  {Proposition}
\newtheorem{lemm+}  [theo+]  {Lemma}
\newtheorem{cor+}  [theo+]  {Corollary}
\newenvironment{theorem}{\begin{theo+}}{\end{theo+}}

\newenvironment{corollary}{\begin{cor+}}{\end{cor+}}

\newcommand{\Ga}{\Gamma}
\newcommand{\De}{\Delta}
\newcommand{\la}{\lambda}

\begin{document}
\baselineskip 18pt
\larger[2]
\title[Reduction formulae for hypergeometric functions]
{Reduction formulae for Karlsson--Minton\\
 type hypergeometric functions}
\author{Hjalmar Rosengren}
\address{Department of Mathematics\\ Chalmers University of Technology and 
G\"oteborg University\\SE-412~96 G\"oteborg, Sweden}
\email{hjalmar@math.chalmers.se}

\begin{abstract}
We prove a master theorem for hypergeometric functions of Karlsson--Minton
type, stating that a very general multilateral $\mathrm{U}(n)$ 
Karlsson--Minton type hypergeometric
series  may be reduced to a finite sum.
This identity contains the Karlsson--Minton summation formula and
many of its known generalizations as special cases, and it  also  implies
several ``Bailey-type'' identities for  $\mathrm{U}(n)$ hypergeometric series,
including  multivariable ${}_{10}W_9$ transformations of
Milne and Newcomb and of  Kajihara.
 Even in the one-variable case our identity
is new, and even in this case its proof depends on the theory of
 multivariable hypergeometric series.
\end{abstract}
\maketitle   

\section{Introduction}

At a first glance, the theory of summation and transformation 
formulas for hypergeometric functions (ordinary and basic or ``$q$'')
may appear as an entangled mess of complicated
formulas involving many parameters.
However, since many  different identities may arise
as special or limit cases of a single formula, it is possible to
structure  the  ``space'' of all such identities by
organizing them into hierarchies. For instance, a large number of the
most useful identities may be understood as a ``Bailey hierarchy'',
originating from Bailey's ${}_{10}W_9$ transformation formula
\cite{ba}, \cite[Equation (III.28)]{gr}. 
  
However, there are identities that do not fit  into the
Bailey hierarchy. An example is the Karlsson--Minton summation
formula \cite{mi,ka}
\begin{equation}\label{km}
{}_{r+2}F_{r+1}\left(\begin{matrix}a,b,c_1+m_1,\dots,c_r+m_r\\
b+1,c_1,\dots,c_r\end{matrix}\,;1\right)=\frac{\Ga(b+1)\Ga(1-a)}{\Ga(1+b-a)}
\prod_{i=1}^r\frac{(c_i-b)_{m_i}}{(c_i)_{m_i}},\end{equation}
where the $m_i$ are non-negative integers and $\operatorname{Re}\,(a+|m|)<1$
 (Minton proved this for $a$ a negative integer and 
Karlsson in general). It belongs to a hierarchy of identities for
hypergeometric series with integral parameter
differences; cf.~\cite{c,gk,g,s,se} for related results.
We will refer to such series as  Karlsson--Minton type hypergeometric
series. (Schlosser \cite{s} prefers the acronym IPD type
series, since the Karlsson--Minton formula may be obtained rather
easily from results known much earlier; cf.~\cite{f,fw}. 
However, in our opinion the
term is perfectly justified, since Minton seems to have been the first to call
attention to this type of series.)

The purpose of this paper is to present 
a master identity for series of Karlsson--Minton type, Theorem \ref{t}.
Not only does it contain a large number of  results
from the papers mentioned above as special cases, 
but it also provides a bridge between the Karlsson--Minton hierarchy 
and the Bailey hierarchy. To find this bridge, and even to state 
our theorem, it is necessary to leave the field of one-variable
series and pass to multivariable series.

The multivariable series that arise are so called $\mathrm{U}(n)$
or $A_n$ hypergeometric series. Series of this type were introduced by
Biedenharn, Holman and Louck \cite{hbl}, motivated by the theory
 of $6j$-symbols of the group $\mathrm{SU}(n)$. 
During the last 25 years, 
the theory of $\mathrm{U}(n)$ series 
(and those connected to other classical groups) 
has been developped extensively by
Gustafson, Milne and many others, and it has  been applied to
 problems in representation theory, number theory and combinatorics.

Experts on multivariable series have  argued that 
certain features of one-variable series
are more easily understood within a multivariable framework; for instance, 
that so called very-well-poised  series should be viewed as
 series in two variables $y_1$, $y_2$ with the summation indices restricted
to a line $y_1+y_2=0$.
The  present paper goes further in this direction and
uses multivariable series
as a tool for studying one-variable series. 

To give an example, a very
degenerate case of Theorem \ref{t} is the following generalization of 
\eqref{km}:
\begin{equation}\label{kmg}\begin{split}&
{}_{r+2}F_{r+1}\left(\begin{matrix}a,b,c_1+m_1,\dots,c_r+m_r\\
d,c_1,\dots,c_r\end{matrix}\,;1\right)
=\frac{\Ga(d)\Ga(d-a-b)}{\Ga(d-a)\Ga(d-b)}
\prod_{i=1}^r\frac{(c_i+1-d)_{m_i}}{(c_i)_{m_i}}\\
&\quad
\times\frac{(a)_{|m|}}{(1+a+b-d)_{|m|}}
\sum_{x_1,\dots,x_r=0}^{m_1,\dots,m_r}\Bigg(\prod_{1\leq i<j\leq r}
\frac{c_i+x_i-c_j-x_j}{c_i-c_j}\frac{(b+1-d)_{|x|}}
{(1-|m|-a)_{|x|}}\\
&\quad\times\prod_{i=1}^r\frac{(c_i-a)_{x_i}}{(1+c_i-d)_{x_i}}
\prod_{i,k=1}^r\frac{(c_i-c_k-m_k)_{x_i}}{(1+c_i-c_k)_{x_i}}\Bigg),
\end{split}\end{equation}
where $\operatorname{Re}\,(a+|m|+b-d)<0$.
When $d=b+1$ all terms in the sum vanish except the one with 
$x_i\equiv 0$, which is $1$, so that we recover \eqref{km};
note also that the case $m_i\equiv 0$ is
 Gauss' classical ${}_2F_1$ summation.  
 The point is that, although \eqref{kmg} is certainly
an interesting identity
from the viewpoint of one-variable series alone,  the finite sum 
on the right is precisely a $\mathrm{U}(n)$ hypergeometric sum, and the 
identity seems difficult to prove without
using the multivariable theory. 

We remark that Karlsson's proof of \eqref{km} 
 gives an alternative expression for
the left-hand side of \eqref{kmg} as a finite 
 sum. Explicitly, one has
\begin{equation}\label{us}\begin{split}&
{}_{r+2}F_{r+1}\left(\begin{matrix}a,b,c_1+m_1,\dots,c_r+m_r\\
d,c_1,\dots,c_r\end{matrix}\,;1\right)
=\frac{\Ga(d)\Ga(d-a-b)}{\Ga(d-a)\Ga(d-b)}\\
&\quad
\times\sum_{x_1,\dots,x_r=0}^{m_1,\dots,m_r}
\frac{(a)_{|x|}(b)_{|x|}}{(1+a+b-d)_{|x|}(c_r)_{|x|}}
\prod_{i=1}^r\frac{(-m_i)_{x_i}}{x_i!}\prod_{i=1}^{r-1}
\frac{(c_{i+1}+m_{i+1})_{x_1+\dots+x_i}}{(c_i)_{x_1+\dots+x_i}}.
\end{split}\end{equation}
This  sum is of a type that is less symmetric than $\mathrm{U}(n)$ series
and probably without much independent interest.

To see how $\mathrm{U}(n)$ series are related to Karlsson--Minton type 
series 
we need only know that the former (we stick to the
classical rather than the $q$-case)
are characterized by the factor 
$$\frac{\De(z+y)}{\De(z)}=\prod_{1\leq i<j\leq n}
\frac{z_i+y_i-z_j-y_j}{z_i-z_j},$$
 where  the $y_i$ are summation indices. 
Suppose we restrict the summation to the line where  
$y_2=\dots=y_{n}=0$ and put $y_1=k$. Then 
$$\frac{\De(z+y)}{\De(z)}=\prod_{i=2}^{n}\frac{z_1+k-z_i}{z_1-z_i}
=\prod_{i=2}^{n}\frac{(z_1-z_i+1)_k}{(z_1-z_i)_k}. $$
If we now choose the parameters $z_i$ so that
$$(z_1-z_2,\dots,z_1-z_{n})=(c_1,c_1+1,\dots,c_1+m_1-1,c_2,
\dots,c_2+m_2-1,\dots,c_r+m_{r}-1), $$
where the $m_i$ are non-negative integers with $|m|=n-1$,  we obtain
$$\frac{\De(z+y)}{\De(z)}=\prod_{i=1}^{r}\frac{(c_i+m_i)_{k}}{(c_i)_k},$$
which is the factor characterizing Karlsson--Minton type series.
So we may view such a series as the restriction of a $\mathrm{U}(n)$ series
to a one-dimensional subspace.

To exploit this observation we must work 
 with $\mathrm{U}(n)$ series for which   restriction of the summation
indices to  lower-dimensional subspaces gives something nice. 
We choose as 
our starting point Gustafson's $\mathrm{U}(n)$ Bailey sum \cite{gu};
cf.~equation \eqref{gi} below. 
In this sum the summation indices live on a hyperplane $y_1+\dots+y_n=0$.
We specialize the parameters so that the terms with $y_n<0$ vanish. 
It then turns out that the sum with $y_n\geq 1$ is of the same type as
the original sum with $y_n\geq 0$. The difference of these two sums
gives the restriction of the original sum
to the space where $y_n=0$, which hence may be computed.
 Iterating this procedure we eventually find an identity
 for the Karlsson--Minton type series obtained 
 by restriction to a one-dimensional subspace. 
Moreover, the previous steps in the iteration give 
identities reducing multivariable Karlsson--Minton series
to finite sums; these are also contained in Theorem \ref{t}.

The details of this derivation are worked out in Section \ref{s1}.
In Section \ref{s2} we state a number of corollaries 
to Theorem \ref{t}. These include summation and transformation formulas
for one- and multivariable Karlsson--Minton type series from 
\cite{c,gk,g,s,se} and also a multivariable generalization of Shukla's
${}_8\psi_8$ summation due to Schlosser \cite{ss}. 
We point out some interesting cases corresponding to lower 
level identities: if  Theorem \ref{t} is on the ${}_6\psi_6$ level,
Corollary \ref{c1p} gives the   ${}_1\psi_1$ version while
Corollaries \ref{2l} and \ref{3l} correspond to the ${}_2\phi_1$ level;
Corollary \ref{2l} is  the $q$-analogue of \eqref{kmg}.
We also indicate how Theorem \ref{t} is related to some  
``Bailey-type'' results for $\mathrm{U}(n)$ series. 
When the Karlsson--Minton type series reduces to a finite sum,
 we may remove the condition that the $m_i$ are non-negative
integers by a polynomial argument,  and recover 
 $\mathrm{U}(n)$ Watson and Bailey transformations, Corollaries \ref{kc}
and \ref{kt},
 which were recently found by Kajihara \cite{k}.
Another interesting case is when the 
Karlsson--Minton type series is one-dimensional. It then has a 
symmetry which is not apparent for the finite sum; this implies
 $\mathrm{U}(n)$ Sears and Bailey transformations due to Milne and Newcomb,
Corollaries \ref{mnc} and \ref{mnb}. 
In Section \ref{q1s} we write down the analogues of Theorem \ref{t}
obtained using instead of 
Gustafson's $\mathrm{U}(n)$ ${}_6\psi_6$ sum the closely related 
$\mathrm{U}(n)$ ${}_5H_5$ and
${}_2H_2$ summation formulas from \cite{gu}.

We finally remark that it seems worthwhile to repeat the analysis
of the present paper starting from summation formulas
different from those  used here.
In fact, in the subsequent paper \cite{r} we apply the $C_n$
${}_6\psi_6$ sum   from \cite{gu2} to obtain a reduction formula
for Karlsson--Minton type hypergeometric series on the root system
$C_n$.

{\bf Acknowledgements:} I would like to thank Michael Schlosser for his
comments on an earlier version of this paper; in particular his
suggestion that I try to derive Theorem 4.6 of \cite{s} using the ideas
of the present paper lead to the inclusion of the more general
Corollary \ref{cmn} below. I also thank Yasushi Kajihara for providing me
with  the manuscript of  \cite{k}.

\section{Notation and a single preliminary}

In the rest of the paper
we will work with  $q$-series, and only discuss the limit
case of classical hypergeometric series briefly in Section \ref{q1s}. 
The base $q$ will be a fixed complex number with 
$0<|q|<1$. We will use the standard notation of \cite{gr}, but since
$q$ is fixed we suppress it from the notation. Thus we write
(note that this is different from the notation  used in the introduction;
cf.~\eqref{sn} below)
\begin{equation}\label{qp}
(a)_k=\begin{cases}(1-a)(1-aq)\dotsm(1-aq^{k-1}),& k\geq 0,\\
\displaystyle\frac{ 1}
{(1-aq^{-1})(1-aq^{-2})\dotsm(1-aq^k)},& k<0,\end{cases}
\end{equation}
$$(a_1,\dots,a_m)_k=(a_1)_k\dotsm(a_m)_k,$$
and analogously for infinite products $(a)_\infty=\prod_{j=0}^\infty(1-aq^j)$.

For $z=(z_1,\dots,z_n)\in\mathbb C^n$  we write $|z|=z_1+\dots+z_n$ and use the
corresponding capital letter to denote the product of the coordinates:
$Z=z_1\dotsm z_n$.

To prove our main theorem all 
we need is Gustafson's multivariable Bailey sum \cite[Theorem 1.15]{gu}, 
which we write as
\begin{equation}\label{gi}
\sum_{\substack{y_1,\dots,y_n=-\infty\\y_1+\dots+y_n=0}}^\infty
\frac{\De(zq^y)}{\De(z)}
\prod_{i,k=1}^n\frac{(a_iz_k)_{y_k}}
{(b_iz_k)_{y_k}}=
\frac{(q/AZ,q^{1-n}BZ)_\infty}{(q,q^{1-n}B/A)_\infty}
\prod_{i,k=1}^n\frac{(b_i/a_k,qz_k/z_i)_\infty}{(q/a_kz_i,b_iz_k)_\infty},
 \end{equation}
where
$$\frac{\De(zq^y)}{\De(z)}=\prod_{1\leq i<j\leq n}
\frac{z_iq^{y_i}-z_jq^{y_j}}{z_i-z_j}.$$
This holds for $|q^{1-n}B/A|<1$, as long as no denominators
vanish. When $n=2$, \eqref{gi} is Bailey's ${}_6\psi_6$ summation 
\cite[Equation (II.33)]{gr}. Gustafson's proof
of \eqref{gi} is based on residue calculus and uses non-trivial 
identities for theta functions.

\section{The  theorem}\label{s1}

Our main result is the following identity. We call it a reduction formula,
since it reduces a very general multilateral Karlsson--Minton type series
to a finite sum.

\begin{theorem}\label{t}
Let  $m_i$ be non-negative integers and
$a_i$, $b_i$, $c_i$,  $z_i$ parameters 
 such that $|q^{1-|m|-n}B/A|<1$ and
none of the denominators in \eqref{te} vanishes.
Then the following identity holds:
\begin{multline}\label{te}
\begin{split}&
\sum_{\substack{y_1,\dots,y_n=-\infty\\y_1+\dots+y_n=0}}^\infty
\frac{\De(zq^y)}{\De(z)}
\prod_{\substack{1\leq k\leq n\\1\leq i\leq p}}\frac{(c_iz_kq^{m_i})_{y_k}}
{(c_iz_k)_{y_k}}\prod_{i,k=1}^n\frac{(a_iz_k)_{y_k}}
{(b_iz_k)_{y_k}}\\
&=\quad
\frac{(q^{1-|m|}/AZ,q^{1-n}BZ)_\infty}{(q,q^{1-|m|-n}B/A)_\infty}
\prod_{i,k=1}^n\frac{(b_i/a_k,qz_k/z_i)_\infty}{(q/a_kz_i,b_iz_k)_\infty}
\prod_{\substack{1\leq k\leq n\\1\leq i\leq p}}\frac{(q^{-m_i}b_k/c_i)_{m_i}}
{(q^{1-m_i}/c_iz_k)_{m_i}}\\
&\quad\quad\times
\sum_{x_1,\dots,x_p=0}^{m_1,\dots,m_p}\frac{\De(cq^x)}{\De(c)}\,q^{|x|}
\frac{(q^{n}/BZ)_{|x|}}{(q^{1-|m|}/AZ)_{|x|}}
\prod_{\substack{1\leq k\leq n\\1\leq i\leq p}}\frac{(c_i/a_k)_{x_i}}
{(qc_i/b_k)_{x_i}}\prod_{i,k=1}^p\frac{(q^{-m_k}c_i/c_k)_{x_i}}
{(qc_i/c_k)_{x_i}}.
\end{split}\end{multline}
\end{theorem}

The condition $|q^{1-|m|-n}B/A|<1$ ensures that the series on the
left-hand side converges absolutely, so that the series manipulations
occurring in the proof are justified. This can be seen exactly as 
in \cite{gu}; we will not discuss questions
of convergence any further.

\begin{proof}
We will prove the following seemingly more 
general identity:
\begin{multline}\label{ten}\begin{split}&
\sum_{\substack{y_1,\dots,y_n=-\infty\\y_1+\dots+y_n=N}}^\infty
\frac{\De(zq^y)}{\De(z)}
\prod_{\substack{1\leq k\leq n\\1\leq i\leq p}}\frac{(c_iz_kq^{m_i})_{y_k}}
{(c_iz_k)_{y_k}}\prod_{i,k=1}^n\frac{(a_iz_k)_{y_k}}
{(b_iz_k)_{y_k}}=q^{\binom{N}{2}}(-q^{|m|}AZ)^N\\
&\quad\times
\frac{(q^{1-|m|-N}/AZ,q^{1+N-n}BZ)_\infty}{(q,q^{1-|m|-n}B/A)_\infty}
\prod_{i,k=1}^n\frac{(b_i/a_k,qz_k/z_i)_\infty}{(q/a_kz_i,b_iz_k)_\infty}
\prod_{\substack{1\leq k\leq n\\1\leq i\leq p}}\frac{(q^{-m_i}b_k/c_i)_{m_i}}
{(q^{1-m_i}/c_iz_k)_{m_i}}\\
&\quad\times
\sum_{x_1,\dots,x_p=0}^{m_1,\dots,m_p}\frac{\De(cq^x)}{\De(c)}\,q^{|x|}
\frac{(q^{n-N}/BZ)_{|x|}}{(q^{1-|m|-N}/AZ)_{|x|}}
\prod_{\substack{1\leq k\leq n\\1\leq i\leq p}}\frac{(c_i/a_k)_{x_i}}
{(qc_i/b_k)_{x_i}}\prod_{i,k=1}^p\frac{(q^{-m_k}c_i/c_k)_{x_i}}
{(qc_i/c_k)_{x_i}}.
\end{split}\end{multline}
Although we do not need it, it is not hard to check that \eqref{ten}
 is  equivalent to \eqref{te}. 
Having the extra parameter $N$ will slightly 
simplify the proof and also be useful later.

We first prove \eqref{ten} for $m_i\equiv 1$ by induction on $p$.
 The starting point $p=0$ is equivalent to \eqref{gi} by a change of 
summation variables;
cf.~Section 5 in \cite{gu}. Then we show that the case
of general $m_i$ may be reduced to the special case $m_i\equiv 1$.

For the first part of the proof we assume that \eqref{ten} holds 
with $n$ replaced by $n+1$ and with $m_1=\dots=m_p=1$. We also specialize
to the case $b_{n+1}=q/z_{n+1}$.
Then the factor $1/(b_{n+1}z_{n+1})_{y_{n+1}}$ on the left-hand side vanishes
unless $y_{n+1}\geq 0$, so that the series is supported on a half-space.
Next we let $a_{n+1}\rightarrow q/z_{n+1}$, which  corresponds
to a removable singularity. After cancelling some factors,  
the first double product on the right-hand side may be 
written as
$$\prod_{i,k=1}^{n+1}\frac{(b_i/a_k,qz_k/z_i)_\infty}
{(q/a_kz_i,b_iz_k)_\infty}=
\prod_{i,k=1}^n\frac{(b_i/a_k,qz_k/z_i)_\infty}{(q/a_kz_i,b_iz_k)_\infty}
\prod_{i=1}^n\frac{1-q^{-1}b_iz_{n+1}}{1-z_{n+1}/z_i},$$
and after further simplifications we obtain the identity
\begin{multline*}\begin{split}
S&=\sum_{\substack{y_1,\dots,y_{n+1}\in\mathbb Z,\ y_{n+1}\geq 0
\\y_1+\dots+y_{n+1}=N}}
\frac{\De(zq^y)}{\De(z)}
\prod_{\substack{1\leq k\leq {n+1}\\1\leq i\leq p}}\frac{(c_iz_kq)_{y_k}}
{(c_iz_k)_{y_k}}\prod_{\substack{1\leq k\leq n+1\\1\leq i\leq n}}
\frac{(a_iz_k)_{y_k}}
{(b_iz_k)_{y_k}}\\
&=q^{\binom N2}(-q^{p+1}\hat A\hat Z)^N
\frac{(q^{-p-N}/\hat A\hat Z,q^{1+N-n}\hat B\hat Z)_\infty}
{(q,q^{-p-n}\hat B/\hat A)_\infty}
\prod_{i,k=1}^n\frac{(b_i/a_k,qz_k/z_i)_\infty}{(q/a_kz_i,b_iz_k)_\infty}\\
&\quad\times
\prod_{i=1}^n\frac{1-q^{-1}b_iz_{n+1}}{1-z_{n+1}/z_i}
\prod_{\substack{1\leq k\leq n\\1\leq i\leq p}}\frac{1-q^{-1}b_k/c_i}
{1-1/c_iz_k}
\sum_{x_1,\dots,x_p=0}^1\Bigg(\frac{\De(cq^x)}{\De(c)}\,q^{|x|}
\frac{(q^{n-N}/\hat B\hat Z)_{|x|}}{(q^{-p-N}/\hat A\hat Z)_{|x|}}\\
&\quad\times
\prod_{\substack{1\leq k\leq n\\1\leq i\leq p}}\frac{(c_i/a_k)_{x_i}}
{(qc_i/b_k)_{x_i}}\prod_{i=1}^p \frac{(q^{-1}c_iz_{n+1})_{x_i}}
{(c_iz_{n+1})_{x_i}}
\prod_{i,k=1}^p\frac{(q^{-1}c_i/c_k)_{x_i}}
{(qc_i/c_k)_{x_i}}\Bigg),
\end{split}\end{multline*}
where $\hat A=a_1\dotsm a_n$ and similarly for $\hat B$ and $\hat Z$.

We now divide the sum into two parts as
$$S=\sum_{y_{n+1}=0}+\sum_{y_{n+1}\geq 1}=S_1+S_2. $$
 By a change of summation variables, $S_2$
may  be reduced to a sum of the same type as $S$.
Indeed, choosing $w_i=z_i$ for $1\leq i\leq n$ and $w_{n+1}=qz_{n+1}$,
so that
$$\frac{\De(w)}{\De(z)}
=\prod_{i=1}^n\frac{1-qz_{n+1}/z_i}{1-z_{n+1}/z_i},$$
we have 
 \begin{multline*}
\begin{split}S_2&
=\sum_{\substack{y_1,\dots,y_{n+1}\in\mathbb Z,\ y_{n+1}\geq 1
\\y_1+\dots+y_{n+1}=N}}
\frac{\De(zq^y)}{\De(z)}
\prod_{\substack{1\leq k\leq {n+1}\\1\leq i\leq p}}\frac{(c_iz_kq)_{y_k}}
{(c_iz_k)_{y_k}}\prod_{\substack{1\leq k\leq n+1\\1\leq i\leq n}}
\frac{(a_iz_k)_{y_k}}
{(b_iz_k)_{y_k}}\\
&=\prod_{i=1}^p\frac{1-c_iz_{n+1}q}{1-c_iz_{n+1}}
\prod_{i=1}^n\left(\frac{1-a_iz_{n+1}}{1-b_iz_{n+1}}\frac{1-qz_{n+1}/z_i}
{1-z_{n+1}/z_i}\right)\\
&\quad\times\sum_{\substack{y_1,\dots,y_{n+1}\in\mathbb Z,\ y_{n+1}\geq 0
\\y_1+\dots+y_{n+1}=N-1}}
\frac{\De(wq^y)}{\De(w)}
\prod_{\substack{1\leq k\leq {n+1}\\1\leq i\leq p}}\frac{(c_iw_kq)_{y_k}}
{(c_iw_k)_{y_k}}\prod_{\substack{1\leq k\leq n+1\\1\leq i\leq n}}
\frac{(a_iw_k)_{y_k}}
{(b_iw_k)_{y_k}}
\\
&=q^{\binom{N-1}{2}}(-q^{p+1}\hat A\hat Z)^{N-1}
\frac{(q^{1-p-N}/\hat A\hat Z,q^{N-n}\hat B\hat Z)_\infty}
{(q,q^{-p-n}\hat B/\hat A)_\infty}
\prod_{i,k=1}^n\frac{(b_i/a_k,qz_k/z_i)_\infty}{(q/a_kz_i,b_iz_k)_\infty}
\\
&\quad\times\prod_{i=1}^n\frac{1-a_iz_{n+1}}{1-z_{n+1}/z_i}
\prod_{\substack{1\leq k\leq n\\1\leq i\leq p}}
\frac{1-q^{-1}b_k/c_i}{1-1/c_iz_k}
\prod_{i=1}^p\frac{1-c_iz_{n+1}q}{1-c_iz_{n+1}}
\sum_{x_1,\dots,x_p=0}^1\Bigg(\frac{\De(cq^x)}{\De(c)}
\\
&\quad\times q^{|x|}
\frac{(q^{n+1-N}/\hat B\hat Z)_{|x|}}{(q^{1-p-N}/\hat A\hat Z)_{|x|}}
\prod_{\substack{1\leq k\leq n\\1\leq i\leq p}}
\frac{(c_i/a_k)_{x_i}}
{(qc_i/b_k)_{x_i}}\prod_{i=1}^p \frac{(c_iz_{n+1})_{x_i}}
{(c_iz_{n+1}q)_{x_i}}
\prod_{i,k=1}^p\frac{(q^{-1}c_i/c_k)_{x_i}}
{(qc_i/c_k)_{x_i}}\Bigg).
\end{split}\end{multline*}
Thus, $S_1$ may be expressed as a finite sum of the form
\begin{equation}\label{se}
S_1=S-S_2=\sum_{x_1,\dots,x_p=0}^1+\sum_{x_1,\dots,x_p=0}^1.\end{equation}

 Next we observe that,
with $y_{n+1}=0$ and $z_{n+1}=1/c_{p+1}$,
\begin{equation*}\begin{split}\frac{\De(zq^y)}{\De(z)}&=\prod_{1\leq i<j\leq n}
\frac{z_iq^{y_i}-z_jq^{y_j}}{z_i-z_j}\prod_{k=1}^n
\frac{z_kq^{y_k}-c_{p+1}^{-1}}{z_k-c_{p+1}^{-1}}\\
&=\prod_{1\leq i<j\leq n}
\frac{z_iq^{y_i}-z_jq^{y_j}}{z_i-z_j}\prod_{k=1}^n
\frac{(c_{p+1}z_kq)_{y_k}}{(c_{p+1}z_k)_{y_k}},\end{split}\end{equation*}
so that
 $$S_1=\sum_{\substack{y_1,\dots,y_n=-\infty\\y_1+\dots+y_n=N}}^\infty
\frac{\De(zq^y)}{\De(z)}
\prod_{\substack{1\leq k\leq n\\1\leq i\leq p+1}}\frac{(c_iz_kq)_{y_k}}
{(c_iz_k)_{y_k}}\prod_{i,k=1}^n\frac{(a_iz_k)_{y_k}}
{(b_iz_k)_{y_k}}, $$ 
 a sum as in \eqref{ten}, still with $m_i\equiv 1$
but with $p$ replaced by $p+1$. 
Writing the corresponding right-hand side of \eqref{ten} as
$$\sum_{x_1,\dots,x_{p+1}=0}^1=\sum_{\substack{0\leq x_1,\dots,x_p\leq 1\\
x_{p+1}=0}}+\sum_{\substack{0\leq x_1,\dots,x_p\leq 1\\
x_{p+1}=1}},$$
it is straight-forward to check that it agrees termwise with \eqref{se}.
Thus, \eqref{ten} holds for $S_1$, and by induction for all $p$ as
long as $m_i\equiv 1$.

To remove the condition $m_i\equiv 1$, we first note that we may
assume $m_i\geq 1$, since if $m_i=0$ all factors involving $m_i$
cancel.
On the other hand, if $m_i> 1$ we may write
$$\frac{(c_iz_kq^{m_i})_{y_k}}{(c_iz_k)_{y_k}}=
\frac{(c_iz_kq)_{y_k}}{(c_iz_k)_{y_k}}
\frac{(c_iz_kq^{2})_{y_k}}{(c_iz_kq)_{y_k}}
\dotsm\frac{(c_iz_kq^{m_i})_{y_k}}{(c_iz_kq^{m_i-1})_{y_k}}, $$
which gives a reduction to the case $m_i\equiv 1$, 
 with $c=(c_1,\dots,c_p)$ replaced by 
\begin{equation}\label{cd}
d=(c_1,qc_1,\dots,q^{m_1-1}c_1,\dots,c_p,qc_p,
\dots,q^{m_p-1}c_p).\end{equation}

We now observe that, when $m_i\equiv 1$, the right-hand side of \eqref{ten} 
contains the factor
\begin{equation}\label{cs}\begin{split}&\frac{\De(cq^x)}{\De(c)}
\prod_{i,k=1}^p\frac{(q^{-1}c_i/c_k)_{x_i}}
{(qc_i/c_k)_{x_i}}\\
&\quad=\prod_{k=1}^n\frac{(q^{-1})_{x_k}}{(q)_{x_k}}
\prod_{1\leq i<j\leq p}\frac{c_iq^{x_i}-c_jq^{x_j}}{c_i-c_j}
\frac{(q^{-1}c_i/c_j)_{x_i}(q^{-1}c_j/c_i)_{x_j}}
{(qc_i/c_j)_{x_i}(qc_j/c_i)_{x_j}}\\
&\quad=(-1)^{|x|}q^{-|x|}\prod_{1\leq i<j\leq p}
\frac{c_iq^{-x_i}-c_jq^{-x_j}}{c_i-c_j},\end{split}\end{equation}
as is easily seen by considering the four cases $x_i,\,x_j=0,\,1$ separately.
If $c$ is replaced by $d$ as above, this expression vanishes unless
$x$ is of the form 
\begin{equation}\label{xy}
x=(\underbrace{\overbrace{1,\dots,1}^{y_1},0,\dots,0}_{m_1},
\dots,\underbrace{\overbrace{1,\dots,1}^{y_p},0,\dots,0}_{m_p}),
\qquad 0\leq y_i\leq m_i. \end{equation}
Rewriting the sum using the $y_i$ as summation variables
and comparing with the corresponding right-hand side of \eqref{ten}, 
we need now only check that, with $c$, $d$ and $x$, $y$ related by 
\eqref{cd} and \eqref{xy}, one has
\begin{multline*}
\prod_{i=1}^{|m|}\frac{1-q^{-1}b_k/d_i}
{1-1/d_iz_k}\frac{\De(dq^x)}{\De(d)}
\prod_{\substack{1\leq k\leq n\\1\leq i\leq |m|}}\frac{(d_i/a_k)_{x_i}}
{(qd_i/b_k)_{x_i}}\prod_{i,k=1}^{|m|}\frac{(q^{-1}d_i/d_k)_{x_i}}
{(qd_i/d_k)_{x_i}}\\
=\prod_{i=1}^{p}\frac{(q^{-m_i}b_k/c_i)_{m_i}}
{(q^{1-m_i}/c_iz_k)_{m_i}}\frac{\De(cq^y)}{\De(c)}
\prod_{\substack{1\leq k\leq n\\1\leq i\leq p}}\frac{(c_i/a_k)_{y_i}}
{(qc_i/b_k)_{y_i}}\prod_{i,k=1}^p\frac{(q^{-m_k}c_i/c_k)_{y_i}}
{(qc_i/c_k)_{y_i}}.
\end{multline*}

Using the obvious identities
$$\prod_{i=1}^{|m|}(ad_i)_{x_i}=\prod_{i=1}^p(ac_i)_{y_i},
\qquad \prod_{i=1}^{|m|}(1-a/d_i)=\prod_{i=1}^p(q^{1-m_i}a/c_i)_{m_i}$$
and \eqref{cs}, we are left with verifying
$$\prod_{1\leq i<j\leq|m|}
\frac{d_iq^{-x_i}-d_jq^{-x_j}}{d_i-d_j}
=(-1)^{|y|}q^{|y|}\frac{\De(cq^y)}{\De(c)}\prod_{i,k=1}^p
\frac{(q^{-m_k}c_i/c_k)_{y_i}}
{(qc_i/c_k)_{y_i}}. $$
To see this we write the left-hand side as
$$\prod_{\{i,\,j:\,i<j,\,x_i=x_j=1\}}q^{-1}
\prod_{\{i,\,j:\,x_i=1,\,x_j=0\}}\frac{d_iq^{-1}-d_j}{d_i-d_j}
=q^{-\binom{|y|}{2}}
\prod_{k,l=1}^p P_{kl}, $$
where
\begin{multline*}\begin{split}P_{kl}&=\prod_{\substack{m_1+\dots+m_{k-1}+1
\leq i\leq m_1+\dots+m_{k-1}+y_k\\m_1+\dots+m_{l-1}+y_l+1
\leq j\leq m_1+\dots+m_{l}}}\frac{d_iq^{-1}-d_j}{d_i-d_j}
=\prod_{\substack{1\leq i\leq y_k\\1\leq j\leq m_{l}-y_l}}
\frac{c_kq^{i-2}-c_lq^{y_l+j-1}}{c_kq^{i-1}-c_lq^{y_l+j-1}}\\
&=\prod_{i=1}^{y_k}\frac{c_kq^{i-1}-c_lq^{m_l}}{c_kq^{i-1}-c_lq^{y_l}}
\,q^{y_l-m_l}=\frac{(q^{-m_l}c_k/c_l)_{y_k}}{(q^{-y_l}c_k/c_l)_{y_k}}.
\end{split}\end{multline*}
It is now enough to show that
$$\prod_{k,l=1}^p\frac{(qc_k/c_l)_{y_k}}{(q^{-y_l}c_k/c_l)_{y_k}}
=(-1)^{|y|}q^{|y|+\binom{|y|}{2}}\prod_{1\leq k<l\leq p}
\frac{c_kq^{y_k}-c_lq^{y_l}}{c_k-c_l}.$$
This follows from elementary identities for $q$-shifted factorials, after
writing
$$ \prod_{k,l=1}^p\frac{(qc_k/c_l)_{y_k}}{(q^{-y_l}c_k/c_l)_{y_k}}
=\prod_{k=1}^p\frac{(q)_{y_k}}{(q^{-y_k})_{y_k}}
\prod_{1\leq k<l\leq p}\frac{(qc_k/c_l)_{y_k}(qc_l/c_k)_{y_l}}
{(q^{-y_l}c_k/c_l)_{y_k}(q^{-y_k}c_l/c_k)_{y_l}}.$$
\end{proof}

\section{Corollaries}\label{s2}

In this section we point out various interesting special cases
and corollaries of Theorem \ref{t}. Throughout it is assumed that
the $m_i$ are non-negative integers and that no denominators in
the identities vanish.

\subsection{Some immediate consequences}

We  first rewrite Theorem \ref{t} in an alternative way, which 
hides much of its symmetry but facilitates comparison with related results
in the literature.
 For this we  replace $n$ by $n+1$ and eliminate
$y_{n+1}=-y_1-\dots-y_n$ from the summation. After the change of variables
\begin{gather*} z_{n+1}\mapsto a^{-1},\qquad
a_{n+1}\mapsto d,\qquad b_{n+1}\mapsto aq/b,\\  
a_i\mapsto e_i/z_i,\qquad b_i\mapsto aq/c_iz_i,\qquad 1\leq i\leq n,\qquad
c_i\mapsto aq/f_i\end{gather*}
we obtain the following  identity.

\begin{corollary}\label{c1}
When $|a^{1+n}q^{1-|m|}/bCdE|<1$, the following identity holds: 
\begin{multline*}
\sum_{y_1,\dots,y_n=-\infty}^\infty
\Bigg(\frac{\De(zq^y)}{\De(z)}\prod_{k=1}^n\frac{1-az_kq^{y_k+|y|}}{1-az_k}
\frac{(b)_{|y|}}{(aq/d)_{|y|}}\prod_{k=1}^n\frac{(c_kz_k)_{|y|}}
{(aqz_k/e_k)_{|y|}}\prod_{i=1}^p\frac{(f_i)_{|y|}}{(q^{-m_i}f_i)_{|y|}}\\
\begin{split}&\quad
\times\prod_{k=1}^n\frac{(dz_k)_{y_k}}{(aqz_k/b)_{y_k}}
\prod_{i,k=1}^n\frac{(e_iz_k/z_i)_{y_k}}{(aqz_k/c_iz_i)_{y_k}}
\prod_{\substack{1\leq k\leq n\\1\leq i\leq p}}
\frac{(aq^{1+m_i}z_k/f_i)_{y_k}}{(aqz_k/f_i)_{y_k}}
\left(\frac{a^{1+n}q^{1-|m|}}{bCdE}\right)^{|y|}\Bigg)\\
&=\frac{(aq^{1-|m|}/dE,a^nq/bC,aq/bd)_\infty}{(a^{n+1}q^{1-|m|}/bCdE,
aq/d,q/b)_\infty}
\prod_{i,k=1}^n\frac{(aqz_k/e_kc_iz_i,qz_k/z_i)_\infty}
{(qz_k/e_kz_i,aqz_k/c_iz_i)_\infty}\\
&\quad\times\prod_{k=1}^n
\frac{(aq/dc_kz_k,q/az_k,aqz_k/be_k,aqz_k)_\infty}
{(q/c_kz_k,q/dz_k,aqz_k/b,aqz_k/e_k)_\infty}
\prod_{\substack{1\leq k\leq n\\1\leq i\leq p}}
\frac{(q^{-m_i}f_i/c_kz_k)_{m_i}}{(q^{-m_i}f_i/az_k)_{m_i}}
\prod_{i=1}^p\frac{(q^{-m_i}f_i/b)_{m_i}}{(q^{-m_i}f_i)_{m_i}}
\\
&\quad\times
\sum_{x_1,\dots,x_p=0}^{m_1,\dots,m_p}\Bigg(
\frac{\De(q^x/f)}{\De(1/f)}\,q^{|x|}
\frac{(bC/a^n)_{|x|}}{(aq^{1-|m|}/dE)_{|x|}}\\
&\quad\times
\prod_{\substack{1\leq k\leq n\\1\leq i\leq p}}\frac{(aqz_k/e_kf_i)_{x_i}}
{(qc_kz_k/f_i)_{x_i}}\prod_{i=1}^p\frac{(aq/df_i)_{x_i}}
{(qb/f_i)_{x_i}}
\prod_{i,k=1}^p\frac{(q^{-m_k}f_k/f_i)_{x_i}}
{(qf_k/f_i)_{x_i}}\Bigg).
\end{split}\end{multline*}
\end{corollary}

The case $p=m_1=1$ of Corollary \ref{c1} is due to Schlosser 
\cite[Theorem 3.4]{ss}. It gives a multivariable analogue of Shukla's
\cite{sh} ${}_8\psi_8$ summation formula, which is obtained if in addition 
$n=1$. Note that in the inductive proof of Theorem \ref{t} given above, the
case $p=m_1=1$ is the first step beyond Gustafson's formula. 
The observation that Schlosser's $n$-variable Shukla summation follows
 from Gustafson's $(n+1)$-variable Bailey summation was in fact the starting
point of the present work. 

In \cite[Theorem 4.2]{s}, 
Schlosser found a transformation formula for series like those
of Corollary \ref{c1}. This identity arises as
 an immediate consequence of Theorem \ref{t},
 namely, from the observation that
the sum on the right-hand side depends only on the product $Z=z_1\dotsm z_n$.
To obtain the 
less compact formulation
in \cite{s} one must first rewrite both series as in Corollary \ref{c1}
 and then replace $y_i$ by $-y_i$ on the right-hand side.

\begin{corollary}[Schlosser]\label{st}
Let $w_i$ and $z_i$ be parameters with
$z_1\dotsm z_n=w_1\dotsm w_n$. Then, assuming $|q^{1-n-|m|}B/A|<1$,
the following identity holds:
 \begin{equation*}\begin{split}&
\sum_{\substack{y_1,\dots,y_n=-\infty\\y_1+\dots+y_n=0}}^\infty
\frac{\De(zq^y)}{\De(z)}
\prod_{\substack{1\leq k\leq n\\1\leq i\leq p}}\frac{(c_iz_kq^{m_i})_{y_k}}
{(c_iz_k)_{y_k}}\prod_{i,k=1}^n\frac{(a_iz_k)_{y_k}}
{(b_iz_k)_{y_k}}\\
&\qquad=
\prod_{i,k=1}^n\frac{(q/a_kw_i,b_iw_k,qz_k/z_i)_\infty}
{(q/a_kz_i,b_iz_k,qw_k/w_i)_\infty}
\prod_{\substack{1\leq k\leq n\\1\leq i\leq p}}\frac{(c_iw_k)_{m_i}}
{(c_iz_k)_{m_i}}\\
&\qquad\quad\times
\sum_{\substack{y_1,\dots,y_n=-\infty\\y_1+\dots+y_n=0}}^\infty
\frac{\De(wq^y)}{\De(w)}
\prod_{\substack{1\leq k\leq n\\1\leq i\leq p}}\frac{(c_iw_kq^{m_i})_{y_k}}
{(c_iw_k)_{y_k}}\prod_{i,k=1}^n\frac{(a_iw_k)_{y_k}}
{(b_iw_k)_{y_k}}
\end{split}\end{equation*}
\end{corollary}

Another immediate corollary of Theorem \ref{t} is obtained
by choosing $BZ=q^n$, which reduces the sum on the right to one term.
Yet again, we obtain an identity of Schlosser \cite[Corollary 4.3]{s}. 

\begin{corollary}[Schlosser]\label{csc}
 If $BZ=q^n$ and $|q^{1-|m|}/AZ|<1$,
 the following identity holds:
 \begin{equation*}\begin{split}&
\sum_{\substack{y_1,\dots,y_n=-\infty\\y_1+\dots+y_n=0}}^\infty
\frac{\De(zq^y)}{\De(z)}
\prod_{\substack{1\leq k\leq n\\1\leq i\leq p}}\frac{(c_iz_kq^{m_i})_{y_k}}
{(c_iz_k)_{y_k}}\prod_{i,k=1}^n\frac{(a_iz_k)_{y_k}}
{(b_iz_k)_{y_k}}\\
&\qquad
=\prod_{i,k=1}^n\frac{(b_i/a_k,qz_k/z_i)_\infty}{(q/a_kz_i,b_iz_k)_\infty}
\prod_{\substack{1\leq k\leq n\\1\leq i\leq p}}\frac{(qc_i/b_k)_{m_i}}
{(c_iz_k)_{m_i}}.
\end{split}\end{equation*}
\end{corollary}

When $n=2$, this is an identity of Chu \cite{c} 
(equivalently, Chu's identity is the case $a=bc$ of Corollary \ref{c2}
below).
In \cite{s}, Corollary \ref{csc} is derived from Corollary \ref{st}
by choosing $w_k=q/b_k$ so that the right-hand side  
is reduced to the term with $y_k\equiv 0$.
We point out that this also happens
under the conditions  $b_k=qa_k$ for $1\leq k\leq n-1$ and the choice
 $w_k=1/a_k$ for $1\leq k\leq n-1$. Alternatively, we may in this situation
use Corollary \ref{mc} below to sum the right-hand side 
of  \eqref{te}. After writing  $a_n=b$ and  $b_n=d$ we obtain 
with either of these two methods the following summation
formula, which is a  $\mathrm{U}(n)$ generalization
of Chu's identity different from Corollary \ref{csc}.

\begin{corollary} For $|q^{-|m|}d/b|<1$, one has
\begin{multline*}
\sum_{\substack{y_1,\dots,y_n=-\infty\\y_1+\dots+y_n=0}}^\infty
\frac{\De(zq^y)}{\De(z)}
\prod_{\substack{1\leq k\leq n\\1\leq i\leq p}}\frac{(c_iz_kq^{m_i})_{y_k}}
{(c_iz_k)_{y_k}}\prod_{\substack{1\leq k\leq n\\1\leq i\leq n-1}}
\frac{(a_iz_k)_{y_k}}
{(qa_iz_k)_{y_k}}\prod_{k=1}^n\frac{(bz_k)_{y_k}}{(dz_k)_{y_k}}\\
\begin{split}&=\frac{(q/AbZ,AdZ)_\infty\prod_{i,k=1}^{n-1}(qa_k/a_i)_\infty
\prod_{i,k=1}^n(qz_k/z_i)_\infty\prod_{k=1}^{n-1}(qa_k/b,d/a_k)_\infty}
{(q)_\infty\prod_{1\leq k\leq n-1,\,1\leq i\leq n}
(q/a_kz_i,qa_kz_i)_\infty
\prod_{k=1}^n(q/bz_k,dz_k)_\infty}\\
&\quad\times\prod_{i=1}^p\frac{(c_iAZ,c_i/a_1,\dots,c_i/a_{n-1})_{m_i}}
{(c_iz_1,\dots,c_iz_n)_{m_i}}.\end{split}
\end{multline*}
\end{corollary}

An important special 
 case of Theorem \ref{t} is when $b_iz_i=q$ for $1\leq i\leq n-1$.
Then only the terms with $y_i\geq 0$   for $1\leq i\leq n-1$ are non-zero,
so that we obtain a multivariable generalization of the unilateral 
$\phi$-series rather than the bilateral $\psi$-series.
When exhibiting this case explicitly we prefer to start from Corollary 
\ref{c1}, where we put $c_1=\dots= c_n=a$.

\begin{corollary}\label{pk}
When $|aq^{1-|m|}/bdE|<1$, the following identity holds: 
\begin{multline*}
\sum_{y_1,\dots,y_n=0}^\infty
\Bigg(\frac{\De(zq^y)}{\De(z)}\prod_{k=1}^n\frac{1-az_kq^{y_k+|y|}}{1-az_k}
\frac{(b)_{|y|}}{(aq/d)_{|y|}}\prod_{k=1}^n\frac{(az_k)_{|y|}}
{(aqz_k/e_k)_{|y|}}\prod_{i=1}^p\frac{(f_i)_{|y|}}{(q^{-m_i}f_i)_{|y|}}
\\
\begin{split}&\quad
\times\prod_{k=1}^n\frac{(dz_k)_{y_k}}{(aqz_k/b)_{y_k}}
\prod_{i,k=1}^n\frac{(e_iz_k/z_i)_{y_k}}{(qz_k/z_i)_{y_k}}
\prod_{\substack{1\leq k\leq n\\1\leq i\leq p}}
\frac{(aq^{1+m_i}z_k/f_i)_{y_k}}{(aqz_k/f_i)_{y_k}}
\left(\frac{aq^{1-|m|}}{bdE}\right)^{|y|}\Bigg)\\
&=\frac{(aq^{1-|m|}/dE,aq/bd)_\infty}{(aq^{1-|m|}/bdE,
aq/d)_\infty}\prod_{k=1}^n
\frac{(aqz_k/be_k,aqz_k)_\infty}
{(aqz_k/b,aqz_k/e_k)_\infty}
\prod_{i=1}^p\frac{(q^{-m_i}f_i/b)_{m_i}}{(q^{-m_i}f_i)_{m_i}}
\\
&\quad\times
\sum_{x_1,\dots,x_p=0}^{m_1,\dots,m_p}\Bigg(
\frac{\De(q^x/f)}{\De(1/f)}\,q^{|x|}
\frac{(b)_{|x|}}{(aq^{1-|m|}/dE)_{|x|}}\\
&\quad\times
\prod_{\substack{1\leq k\leq n\\1\leq i\leq p}}\frac{(aqz_k/e_kf_i)_{x_i}}
{(aqz_k/f_i)_{x_i}}\prod_{i=1}^p\frac{(aq/df_i)_{x_i}}
{(qb/f_i)_{x_i}}
\prod_{i,k=1}^p\frac{(q^{-m_k}f_k/f_i)_{x_i}}
{(qf_k/f_i)_{x_i}}\Bigg).
\end{split}\end{multline*}
\end{corollary}

Yet another interesting case of Theorem \ref{t} is when the
sum on the right-hand side is supported on a hyperplane $|x|=N$. 
In this case we prefer to start from \eqref{ten} and assume that
 $AZ=q^{-|m|}$, $BZ=q^n$. Then the factor 
$$\frac{(q^{1-|m|-N}/AZ)_\infty(q^{n-N}/BZ)_{|x|}}{(q^{1-|m|-N}/AZ)_{|x|}}
=(q^{1-N+|x|})_\infty(q^{-N})_{|x|}$$
on the right vanishes unless $|x|=N$, which reduces
the finite sum  to a $(p-1)$-variable ${}_{2n+6}W_{2n+5}$
rather than a $p$-variable ${}_{n+2}\phi_{n+1}$.
Moreover, the condition for convergence of the left-hand side reduces
to $|q|<1$, which is automatically satisfied. This leads to the
following identity, which may be viewed as a well-poised version 
of Theorem~\ref{t}.

\begin{corollary}\label{pkb}
If $AZ=q^{-|m|}$ and $BZ=q^n$
the following identity holds:
\begin{multline*}
\sum_{\substack{y_1,\dots,y_n=-\infty\\y_1+\dots+y_n=N}}^\infty
\frac{\De(zq^y)}{\De(z)}
\prod_{\substack{1\leq k\leq n\\1\leq i\leq p}}\frac{(c_iz_kq^{m_i})_{y_k}}
{(c_iz_k)_{y_k}}\prod_{i,k=1}^n\frac{(a_iz_k)_{y_k}}
{(b_iz_k)_{y_k}}
=\prod_{i,k=1}^n\frac{(b_i/a_k,qz_k/z_i)_\infty}{(q/a_kz_i,b_iz_k)_\infty}\\
\times
\prod_{\substack{1\leq k\leq n\\1\leq i\leq p}}\frac{(qc_i/b_k)_{m_i}}
{(c_iz_k)_{m_i}}\sum_{\substack{x_1,\dots,x_p=0\\x_1+\dots+x_p=N}}^
{m_1,\dots,m_p}\frac{\De(cq^x)}{\De(c)}
\prod_{\substack{1\leq k\leq n\\1\leq i\leq p}}\frac{(c_i/a_k)_{x_i}}
{(qc_i/b_k)_{x_i}}\prod_{i,k=1}^p\frac{(q^{-m_k}c_i/c_k)_{x_i}}
{(qc_i/c_k)_{x_i}}.
\end{multline*}
\end{corollary} 

\subsection{Kajihara's transformations}

Theorem \ref{t} may be viewed as a transformation formula relating
hypergeometric series of \emph{different} dimension. Although such results
are rare, some identities of this type
were recently obtained by Kajihara (cf.~\cite{gkr,kr} for other transformations
with this property). 
In fact, it is possible to obtain  Kajihara's identities from Theorem \ref{t}
 by choosing the parameters so that the left-hand side is a finite
sum and then applying a standard ``polynomial argument''.
Starting from the case $b=q^{-N}$ of
Corollary \ref{pk}  we obtain in this way Proposition 6.1
of \cite{k}.
 In the one-variable case
 $n=p=1$ this is a Watson-type transformation that
 may be obtained by combining Equations (III.15) and 
(III.18) of \cite{gr}.

\begin{corollary}[Kajihara]\label{kc}
The following identity holds: 
\begin{multline*}
\sum_{\substack{y_1,\dots,y_n\geq 0\\y_1+\dots+y_n\leq N}}
\Bigg(\frac{\De(zq^y)}{\De(z)}\prod_{k=1}^n\frac{1-az_kq^{y_k+|y|}}{1-az_k}
\frac{(q^{-N})_{|y|}}{(aq/d)_{|y|}}\prod_{k=1}^n\frac{(az_k)_{|y|}}
{(aqz_k/e_k)_{|y|}}\prod_{i=1}^p\frac{(f_i)_{|y|}}{(aq/g_i)_{|y|}}
\\
\begin{split}&\quad
\times\prod_{k=1}^n\frac{(dz_k)_{y_k}}{(aq^{1+N}z_k)_{y_k}}
\prod_{i,k=1}^n\frac{(e_iz_k/z_i)_{y_k}}{(qz_k/z_i)_{y_k}}
\prod_{\substack{1\leq k\leq n\\1\leq i\leq p}}
\frac{(g_iz_k)_{y_k}}{(aqz_k/f_i)_{y_k}}
\left(\frac{a^{p+1}q^{p+N+1}}{dEFG}\right)^{|y|}\Bigg)\\
&=\frac{(a^{p+1}q^{p+1}/dEFG)_N}{(aq/d)_N}\prod_{k=1}^n
\frac{(aqz_k)_N}{(aqz_k/e_k)_N}
\prod_{i=1}^p\frac{(f_i)_{N}}{(aq/g_i)_{N}}\\
&\quad\times
\sum_{\substack{x_1,\dots,x_p\geq 0\\x_1+\dots+x_p\leq N}}\Bigg(
\frac{\De(q^x/f)}{\De(1/f)}\,q^{|x|}
\frac{(q^{-N})_{|x|}}{(a^{p+1}q^{p+1}/dEFG)_{|x|}}\\
&\quad\times
\prod_{\substack{1\leq k\leq n\\1\leq i\leq p}}\frac{(aqz_k/f_ie_k)_{x_i}}
{(aqz_k/f_i)_{x_i}}\prod_{i=1}^p\frac{(aq/df_i)_{x_i}}
{(q^{1-N}/f_i)_{x_i}}
\prod_{i,k=1}^p\frac{(aq/g_kf_i)_{x_i}}
{(qf_k/f_i)_{x_i}}\Bigg).
\end{split}\end{multline*}
\end{corollary} 

\begin{proof}
When $g_i=aq^{1+m_i}/f_i$, this is the case $b=q^{-N}$ of Corollary \ref{pk}.
Initially we have this only when $|q^{|m|}|>|aq^{1+N}/dE|$.
However, since both sides of the identity depend rationally on $d$, 
it extends to all non-negative
integers $m_i$. Thus, if we
 multiply  the identity with $(aq/g_1,\dots,aq/g_p)_N$
both sides will be polynomials in the variables $1/g_i$ that agree
at an infinite number of points and are thus identical.
\end{proof}

If we repeat the same procedure starting from 
 Corollary \ref{pkb}, we recover Kajihara's beautiful
Bailey-type transformation formula \cite[Proposition 6.2]{k}. 
Namely, in the case $b_i=q/z_i$  the left-hand side
of Corollary \ref{pkb} becomes a finite sum. We may then, as 
in the proof of Corollary \ref{kc}, remove the condition that
the $m_i$ are non-negative integers and obtain the following identity. 
When $n=p=2$ it is a ${}_{10}W_9$ transformation
 which may be obtained by iterating \cite[Equation (III.28)]{gr}.
From the proof given above it is apparent that
Corollary \ref{kt} is actually the special case $d=aq^{N}$, 
$EFG=a^pq^{p}$ of Corollary \ref{kc}, but since this is not
noted in \cite{k} we list it as a separate corollary.

\begin{corollary}[Kajihara]\label{kt}
For $W=ABZ$, the following identity holds:
\begin{equation*}\begin{split}&
\sum_{\substack{y_1,\dots,y_n\geq 0\\y_1+\dots+y_n=N}}
\frac{\De(zq^y)}{\De(z)}\prod_{i,k=1}^n\frac{(a_iz_k)_{y_k}}
{(qz_k/z_i)_{y_k}}
\prod_{\substack{1\leq k\leq n\\1\leq i\leq p}}\frac{(b_iz_k)_{y_k}}
{(w_iz_k)_{y_k}}\\
&\qquad=\sum_{\substack{x_1,\dots,x_p\geq 0\\x_1+\dots+x_p=N}}
\frac{\De(wq^x)}{\De(w)}\prod_{i,k=1}^p\frac{(w_i/b_k)_{x_i}}
{(qw_i/w_k)_{x_i}}
\prod_{\substack{1\leq k\leq n\\1\leq i\leq p}}\frac{(w_i/a_k)_{x_i}}
{(w_iz_k)_{x_i}}.\end{split}\end{equation*}
\end{corollary}

\subsection{Low values of $n$}

The right-hand side of \eqref{te} is a $p$-variable analogue
of a balanced ${}_{n+2}\phi_{n+1}$-series. In view of the classical
${}_3\phi_2$ summation and ${}_4\phi_3$ transformation formulas, 
one would expect the cases $n=1$ and $n=2$ to be of special interest.
Indeed, when $n=1$ the left-hand side of Theorem \ref{t} 
 reduces to $1$, so that the multivariable ${}_3\phi_2$  on the
right-hand side can be summed. After replacing $p$ by $n$ and 
relabelling the parameters,
we recover the following multivariable 
$q$-Saalsch\"utz summation formula due to Milne \cite[Theorem 4.1]{m}.
The resulting new proof of Milne's identity is not very natural, 
but it gives a first illustration of
 how Theorem \ref{t} forms a bridge between  different types of
identities for hypergeometric functions. 

\begin{corollary}[Milne]\label{mc}
If $q^{1-|m|}ab=cd$, the following identity holds:
\begin{equation}\label{mce}\begin{split}&
\sum_{x_1,\dots,x_n=0}^{m_1,\dots,m_n}\frac{\De(zq^x)}{\De(z)}\,q^{|x|}
\frac{(a)_{|x|}}{(c)_{|x|}}
\prod_{i=1}^n\frac{(bz_i)_{x_i}}
{(dz_i)_{x_i}}\prod_{i,k=1}^n\frac{(q^{-m_k}z_i/z_k)_{x_i}}
{(qz_i/z_k)_{x_i}}\\
&\qquad=\frac{(d/b)_{|m|}}{(d/ab)_{|m|}}\prod_{i=1}^n
\frac{(dz_i/a)_{m_i}}{(dz_i)_{m_i}}.
\end{split}\end{equation}
\end{corollary}

If we put $n=1$  in Corollary \ref{pkb}, or equivalently multiply
\eqref{mce} with $(d/ab)_{|m|}$ and then let $a=q^{-N}$, $c=q^{1-N}$,
we obtain 
$$ \sum_{\substack{x_1,\dots,x_p=0\\x_1+\dots+x_p=N}}^
{m_1,\dots,m_p}\frac{\De(zq^x)}{\De(z)}
\prod_{i=1}^p\frac{(q^{|m|}dz_i)_{x_i}}
{(dz_i)_{x_i}}\prod_{i,k=1}^p\frac{(q^{-m_k}z_i/z_k)_{x_i}}
{(qz_i/z_k)_{x_i}}=\frac{(q^{-|m|})_N}{(q)_N}\prod_{i=1}^p
\frac{(q^{m_i}dz_i)_{N}}{(dz_i)_N}.$$
By a polynomial argument (cf.~the proof of 
 Corollary \ref{kc}), this is equivalent to a multivariable $q$-Dougall
summation of Milne \cite[Theorem 6.17]{mr} 
(stated in more transparent notation in
\cite[Theorem A.5]{mn}). The observation that Milne's $q$-Dougall
sum follows easily from his $q$-Saalsch\"utz sum appears to be new.

\begin{corollary}[Milne]\label{cmd}
The following identity holds:
$$ \sum_{\substack{x_1,\dots,x_n\geq 0\\x_1+\dots+x_n=N}}
\frac{\De(zq^x)}{\De(z)}
\prod_{i=1}^n\frac{(dz_i/E)_{x_i}}
{(dz_i)_{x_i}}\prod_{i,k=1}^n\frac{(e_kz_i/z_k)_{x_i}}
{(qz_i/z_k)_{x_i}}=\frac{(E)_N}{(q)_N}\prod_{i=1}^n
\frac{(dz_i/e_i)_{N}}{(dz_i)_N}.$$
\end{corollary}

Next we turn to the case
 $n=2$ of Theorem \ref{t}, when the  left-hand side is a 
one-variable very-well-poised 
${}_{p+6}\psi_{p+6}$ series, and the sum on the right 
a $p$-variable terminating balanced ${}_4\phi_3$.
We  write it out  explicitly by
 letting $n=1$ and (without loss of generality) $z_1=1$ 
in  Corollary \ref{c1}.

\begin{corollary}\label{c2}
For $|a^2q^{1-|m|}/bcde|<1$, the following identity holds:
\begin{multline*}\sum_{y=-\infty}^\infty
\frac{1-aq^{2y}}{1-a}\frac{(b,c,d,e)_y}{(aq/b,aq/c,aq/d,aq/e)_y}
\prod_{i=1}^p\frac{(f_i,aq^{1+m_i}/f_i)_{y}}{(q^{-m_i}f_i,aq/f_i)_{y}}
\left(\frac{a^{2}q^{1-|m|}}{bcde}\right)^{y}\\
\begin{split}
&=\frac{(q,aq,q/a,aq/bc,aq/bd,aq/be,aq/cd,aq/ce,aq^{1-|m|}/de)_\infty}
{(q/b,q/c,q/d,q/e,aq/b,aq/c,aq/d,aq/e,a^{2}q^{1-|m|}/bcde)_\infty}\\
&\quad\times\prod_{i=1}^p\frac{(q^{-m_i}f_i/b,q^{-m_i}f_i/c)_{m_i}}
{(q^{-m_i}f_i,q^{-m_i}f_i/a)_{m_i}}
\sum_{x_1,\dots,x_p=0}^{m_1,\dots,m_p}\Bigg(
\frac{\De(q^x/f)}{\De(1/f)}\,q^{|x|}
\frac{(bc/a)_{|x|}}{(aq^{1-|m|}/de)_{|x|}}\\
&\quad\times\prod_{i=1}^p\frac{(aq/df_i,aq/ef_i)_{x_i}}
{(qb/f_i,qc/f_i)_{x_i}}\prod_{i,k=1}^p\frac{(q^{-m_k}f_k/f_i)_{x_i}}
{(qf_k/f_i)_{x_i}}\Bigg).
\end{split}\end{multline*}
\end{corollary}

In our opinion,  Corollary \ref{c2} is an interesting
 result both from the viewpoint
of one- and multivariable $q$-series. It is interesting both that the
 bilateral series on the left may be reduced to a finite sum
and that the  multivariable  ${}_4\phi_3$ on the  right
 may be written as a single series. Corollary \ref{c2}
may be compared with Theorem
1.7 of \cite{mo}, which also reduces the left-hand side
  to a finite sum but in a  less symmetric way, similar to \eqref{us}.

If we put $e=a$ in Corollary \ref{c2} or $n=1$ in Corollary \ref{pk},
we obtain the following unilateral identity.

\begin{corollary}\label{c3}
For $|aq^{1-|m|}/bcd|<1$, the following identity holds:
\begin{multline*}
\sum_{y=0}^\infty
\frac{1-aq^{2y}}{1-a}\frac{(a,b,c,d)_y}{(q,aq/b,aq/c,aq/d)_y}
\prod_{i=1}^p\frac{(f_i,aq^{1+m_i}/f_i)_{y}}{(q^{-m_i}f_i,aq/f_i)_{y}}
\left(\frac{aq^{1-|m|}}{bcd}\right)^{y}\\
\begin{split}&\quad=\frac{(aq,aq/bc,aq/bd,aq/cd)_\infty}
{(aq/b,aq/c,aq/d,aq^{1-|m|}/bcd)_\infty}
\prod_{i=1}^p\frac{(q^{-m_i}f_i/b,q^{-m_i}f_i/c)_{m_i}}
{(q^{-m_i}f_i,q^{-m_i}f_i/a)_{m_i}}\,(q^{1-|m|}/d)_{|m|}\\
&\quad\quad\times\sum_{x_1,\dots,x_p=0}^{m_1,\dots,m_p}
\frac{\De(q^x/f)}{\De(1/f)}\,q^{|x|}
\frac{(bc/a)_{|x|}}{(q^{1-|m|}/d)_{|x|}}
\prod_{i=1}^p\frac{(aq/df_i,q/f_i)_{x_i}}
{(qb/f_i,qc/f_i)_{x_i}}\prod_{i,k=1}^p\frac{(q^{-m_k}f_k/f_i)_{x_i}}
{(qf_k/f_i)_{x_i}}.
\end{split}\end{multline*}
\end{corollary}

The case $p=1$ of Corollary \ref{c3} gives a version of Watson's 
transformation formula,
relating a non-terminating ${}_8W_7$ and a terminating ${}_4\phi_3$, 
which may be obtained by combining  equations (III.15)
and (III.20) from \cite{gr}. The case
 $a=bc$   gives an identity of Gasper \cite{g}. Another
interesting case arises if we let $a$, $b$ and $f_i$ tend to $0$ in such a way
that $aq/b$ and $aq/f_i$ are fixed. After  relabelling  the parameters
we obtain the following identity. When $d=bq$, it 
reduces to Gasper's \cite{gk}, \cite[Equation (II.26)]{gr}  
$q$-analogue of 
the Karlsson--Minton summation formula. 

\begin{corollary}\label{2l}
 For $|dq^{-|m|}/ab|<1$, the following identity holds:
\begin{multline*}\sum_{y=0}^\infty
\frac{(a,b)_y}{(q,d)_y}
\prod_{i=1}^p\frac{(c_iq^{m_i})_{y}}{(c_i)_{y}}
\left(\frac{dq^{-|m|}}{ab}\right)^{y}\\
\begin{split}&=\frac{(d/a,d/b)_\infty}
{(d,q^{-|m|}d/ab)_\infty}
\prod_{i=1}^p\frac{(qc_i/d)_{m_i}}
{(c_i)_{m_i}}\left(\frac d q\right)^{|m|}(q^{1-|m|}/a)_{|m|}\\
&\quad\times\sum_{x_1,\dots,x_p=0}^{m_1,\dots,m_p}
\frac{\De(cq^x)}{\De(c)}\,\left(\frac q b\right)^{|x|}
\frac{(qb/d)_{|x|}}{(q^{1-|m|}/a)_{|x|}}
\prod_{i=1}^p\frac{(c_i/a)_{x_i}}
{(qc_i/d)_{x_i}}\prod_{i,k=1}^p\frac{(q^{-m_k}c_i/c_k)_{x_i}}
{(qc_i/c_k)_{x_i}}.
\end{split}\end{multline*}
\end{corollary}

If we put $c=a$ in Corollary \ref{c2} and then let $a$, $b$ and $f_i$ 
tend to  $0$ in such a way that $aq/b$ and $aq/f_i$ are fixed,
we obtain the following identity. 

\begin{corollary}\label{3l}
 For $|dq^{-|m|}/ab|<1$, the following identity holds:
\begin{equation*}\begin{split}&\sum_{y=0}^\infty
\frac{(a,b)_y}{(q,d)_y}
\prod_{i=1}^p\frac{(c_iq^{m_i})_{y}}{(c_i)_{y}}
\left(\frac{dq^{-|m|}}{ab}\right)^{y}=\frac{(d/a,d/b)_\infty}
{(d,q^{-|m|}d/ab)_\infty}
\prod_{i=1}^p(q^{-m_i}d/c_i)_{m_i}\\
&\qquad\times\sum_{x_1,\dots,x_p=0}^{m_1,\dots,m_p}
\frac{\De(cq^x)}{\De(c)}\,q^{|x|}
\prod_{i=1}^p\frac{(c_i/a,c_i/b)_{x_i}}
{(c_i,qc_i/d)_{x_i}}\prod_{i,k=1}^p\frac{(q^{-m_k}c_i/c_k)_{x_i}}
{(qc_i/c_k)_{x_i}}.
\end{split}\end{equation*}
\end{corollary}

 The equality of the left-hand sides of Corollaries \ref{2l} and \ref{3l}
implies a transformation formula between the sums on the right.
More generally, we may start from the observation that
the left-hand side of  Corollary \ref{c2} is invariant under 
interchanging $b$ and $d$. This symmetry is peculiar to the case $n=2$
of Theorem \ref{t}. 
Replacing $p$ by $n$ and relabelling the parameters we obtain the
following  multivariable analogue
of Sears' transformation formula. As is explained below, it is a
degenerate case  of a multivariable ${}_{10}W_9$ transformation
due to Milne and Newcomb. If we let $a$ and $d$ tend to $0$
in Corollary \ref{mnc} with $a/d$  fixed
we recover the $\mathrm{U}(n)$ ${}_3\phi_2$ transformation that 
connects Corollaries \ref{2l} and \ref{3l}.

\begin{corollary}\label{mnc}
If $q^{1-|m|}abc=def$, the following identity holds:
\begin{multline}\label{mnci}
\sum_{x_1,\dots,x_n=0}^{m_1,\dots,m_n}\frac{\De(zq^x)}{\De(z)}\,q^{|x|}
\,\frac{(a)_{|x|}}{(d)_{|x|}}
\prod_{i=1}^n\frac{(bz_i,cz_i)_{x_i}}
{(ez_i,fz_i)_{x_i}}\prod_{i,k=1}^n\frac{(q^{-m_k}z_i/z_k)_{x_i}}
{(qz_i/z_k)_{x_i}}\\
\begin{split}&\quad=\frac{(f/b)_{|m|}}{(q^{1-|m|}/d)_{|m|}}\prod_{i=1}^n
\frac{(q^{1-|m|}bz_i/d)_{m_i}}{(fz_i)_{m_i}}\\
&\quad\quad\times
\sum_{x_1,\dots,x_n=0}^{m_1,\dots,m_n}\frac{\De(zq^x)}{\De(z)}\,q^{|x|}\,
\frac{(e/c)_{|x|}}{(q^{1-|m|}b/f)_{|x|}}
\prod_{i=1}^n\frac{(bz_i,ez_i/a)_{x_i}}
{(ez_i,q^{1-|m|}bz_i/d)_{x_i}}\prod_{i,k=1}^n\frac{(q^{-m_k}z_i/z_k)_{x_i}}
{(qz_i/z_k)_{x_i}}.
\end{split}\end{multline}
\end{corollary}

Finally we consider the ``well-poised'' specialization of 
 Corollary \ref{mnc} (similar to Corollaries \ref{pkb}, \ref{kt} and \ref{cmd}
above). If we multiply \eqref{mnci} with $(q^{1-|m|}/d)_{|m|}$ and then let
$a=q^{-N}$, $d=q^{1-N}$, $c=eq^L$, $f=bq^{L-|m|}$ with $N$ and $L$
non-negative integers, we obtain
\begin{multline}\label{pbt}
\sum_{\substack{x_1,\dots,x_n=0\\x_1+\dots+x_n=N}}^{m_1,\dots,m_n}
\frac{\De(zq^x)}{\De(z)}
\prod_{i=1}^n\frac{(bz_i,q^Lez_i)_{x_i}}
{(ez_i,q^{L-|m|}bz_i)_{x_i}}
\prod_{i,k=1}^n\frac{(q^{-m_k}z_i/z_k)_{x_i}}
{(qz_i/z_k)_{x_i}}=\frac{(q)_{L}(q^{-|m|})_N}{(q)_{N}(q^{-|m|})_L}\\
\times\prod_{i=1}^n
\frac{(q^{N-|m|}bz_i)_{m_i}}{(q^{L-|m|}bz_i)_{m_i}}
\sum_{\substack{x_1,\dots,x_n=0\\x_1+\dots+x_n=L}}^{m_1,\dots,m_n}
\frac{\De(zq^x)}{\De(z)}
\prod_{i=1}^n\frac{(bz_i,q^Nez_i)_{x_i}}
{(ez_i,q^{N-|m|}bz_i)_{x_i}}\prod_{i,k=1}^n\frac{(q^{-m_k}z_i/z_k)_{x_i}}
{(qz_i/z_k)_{x_i}}.
\end{multline}
This corresponds to writing the series
$$\sum_{y=-\infty}^\infty
\frac{1-aq^{2y}}{1-a}\frac{(e,aq^{-N}/e,q^{N+L-|m|}e,aq^{-L}/e)_y}
{(aq/e,q^{N+1}e,aq^{1+|m|-L-N}/e,q^{L+1}e)_y}
\prod_{i=1}^p\frac{(f_i,aq^{1+m_i}/f_i)_{y}}{(q^{-m_i}f_i,aq/f_i)_{y}}
\,q^{y}$$
(with $a=q^{N+L-|m|}be$, $f_i=q/z_i$)
as a finite sum in  two ways using Corollary \ref{c2}.

Equation \eqref{pbt} is a multivariable ${}_{10}W_9$ transformation
closely related to those of Milne and Newcomb \cite{mn}. 
We indicate how to recover \cite[Theorem 3.1]{mn} from \eqref{pbt}.
One should then replace $n$ by $n+1$ in \eqref{pbt} and eliminate
$x_{n+1}$ from the summations. By a polynomial argument, as in 
the proof of  Corollary \ref{kc}, one may first remove the condition 
$m_{n+1}\in\mathbb N$ and then the conditions
$N\in\mathbb N$ and $L\in\mathbb N$. After relabelling the parameters, one 
obtains the following identity.

\begin{corollary}[Milne and Newcomb]\label{mnb}
Assuming that $bcdefg=a^3q^{2+|m|}$ and writing
$\lambda=qa^2/bef$, one has the identity
\begin{multline*}
\sum_{x_1,\dots,x_n=0}^{m_1,\dots,m_n}\Bigg(
\frac{\De(zq^x)}{\De(z)}\,q^{|x|}\prod_{i=1}^n\frac{1-az_iq^{x_i+|x|}}{1-az_i}
\prod_{i=1}^n\frac{(bz_i,cz_i,dz_i)_{x_i}}
{(aqz_i/e,aqz_i/f,aqz_i/g)_{x_i}}\\
\begin{split}&\quad\times\prod_{i,k=1}^n\frac{(q^{-m_k}z_i/z_k)_{x_i}}
{(qz_i/z_k)_{x_i}}\frac{(e,f,g)_{|x|}}{(aq/b,aq/c,aq/d)_{|x|}}\prod_{i=1}^n
\frac{(az_i)_{|x|}}{(q^{1+m_i}az_i)_{|x|}}\Bigg)\\
&=\left(\frac a\lambda\right)^{|m|}\frac{(\la q/c,\la q/d)_{|m|}}
{(aq/c,aq/d)_{|m|}}\prod_{i=1}^n\frac{(aqz_i,\la qz_i/g)_{m_i}}
{(\la qz_i,aqz_i/g)_{m_i}}\\
&\quad\times\sum_{x_1,\dots,x_n=0}^{m_1,\dots,m_n}\Bigg(
\frac{\De(zq^x)}{\De(z)}\,q^{|x|}\prod_{i=1}^n\frac{1-\la 
z_iq^{x_i+|x|}}{1-\la z_i}
\prod_{i=1}^n\frac{(aqz_i/ef,cz_i,dz_i)_{x_i}}
{(aqz_i/e,aqz_i/f,\la qz_i/g)_{x_i}}\\
&\quad\times\prod_{i,k=1}^n\frac{(q^{-m_k}z_i/z_k)_{x_i}}
{(qz_i/z_k)_{x_i}}\frac{(aq/be,aq/bf,g)_{|x|}}{(aq/b,\la q/c,\la q/d)_{|x|}}
\prod_{i=1}^n
\frac{(\la z_i)_{|x|}}{(q^{1+m_i}\la z_i)_{|x|}}\Bigg).
\end{split}\end{multline*}
\end{corollary}

If we let $a$, $d$, $f$, $g\rightarrow 0$ in Corollary \ref{mnb}
in such a way that $a/d$, $a/f$ and $a/g$ are fixed, we recover 
Corollary \ref{mnc}.

\subsection{Generalized ${}_1\psi_1$ series}

Next we consider identities that may be obtained from Theorem \ref{t}
 by multiplying both sides of \eqref{ten} with some function $f(N)$ and then 
summing over $N$. In general, this gives
\begin{multline*}
\sum_{y_1,\dots,y_n=-\infty}^\infty
\frac{\De(zq^y)}{\De(z)}
\prod_{\substack{1\leq k\leq n\\1\leq i\leq p}}\frac{(c_iz_kq^{m_i})_{y_k}}
{(c_iz_k)_{y_k}}\prod_{i,k=1}^n\frac{(a_iz_k)_{y_k}}
{(b_iz_k)_{y_k}}\,f(|y|)\\
\begin{split}&
=\frac{(q^{1-|m|}/AZ,q^{1-n}BZ)_\infty}{(q,q^{1-|m|-n}B/A)_\infty}
\prod_{i,k=1}^n\frac{(b_i/a_k,qz_k/z_i)_\infty}{(q/a_kz_i,b_iz_k)_\infty}
\prod_{\substack{1\leq k\leq n\\1\leq i\leq p}}\frac{(q^{-m_i}b_k/c_i)_{m_i}}
{(q^{1-m_i}/c_iz_k)_{m_i}}\\
&\quad\times
\sum_{x_1,\dots,x_p=0}^{m_1,\dots,m_p}
\Bigg(\frac{\De(cq^x)}{\De(c)}\,q^{|x|}
\frac{(q^n/BZ)_{|x|}}{(q^{1-|m|}/AZ)_{|x|}}
\prod_{\substack{1\leq k\leq n\\1\leq i\leq p}}\frac{(c_i/a_k)_{x_i}}
{(qc_i/b_k)_{x_i}}\prod_{i,k=1}^p\frac{(q^{-m_k}c_i/c_k)_{x_i}}
{(qc_i/c_k)_{x_i}}\\
&\quad\times\sum_{N=-\infty}^\infty \frac{(q^{|m|-|x|}AZ)_N}
{(q^{1-n-|x|}BZ)_N}\,f(N)\Bigg),
\end{split}\end{multline*}
which is valid whenever the series involved are absolutely convergent.
This is of course most interesting when the sum in $N$
simplifies. As an example,
when $f(N)=t^N$ an application of Ramanujan's ${}_1\psi_1$
sum \cite[Equation (II.29)]{gr} yields the following identity, 
which reduces to Gustafson's multivariable ${}_1\psi_1$ sum
\cite[Theorem 1.17]{gu} when $p=0$. 

\begin{corollary}\label{c1p}
For $|q^{1-|m|-n}B/A|<|t|<1$, the following identity holds:
\begin{multline*}
\sum_{y_1,\dots,y_n=-\infty}^\infty
\frac{\De(zq^y)}{\De(z)}
\prod_{\substack{1\leq k\leq n\\1\leq i\leq p}}\frac{(c_iz_kq^{m_i})_{y_k}}
{(c_iz_k)_{y_k}}\prod_{i,k=1}^n\frac{(a_iz_k)_{y_k}}
{(b_iz_k)_{y_k}}\,t^{|y|}\\
\begin{split}&
=\frac{(AZt,q/AZt)_\infty}{(t,q^{1-n}B/At)_\infty}\frac{1}{(q^nAt/B)_{|m|}}
\prod_{i,k=1}^n\frac{(b_i/a_k,qz_k/z_i)_\infty}{(q/a_kz_i,b_iz_k)_\infty}
\prod_{\substack{1\leq k\leq n\\1\leq i\leq p}}\frac{(qc_i/b_k)_{m_i}}
{(c_iz_k)_{m_i}}\\
&\quad\times
\sum_{x_1,\dots,x_p=0}^{m_1,\dots,m_p}\frac{\De(cq^x)}{\De(c)}\,
\left(q^{n+|m|}\frac{At}{B}\right)^{|x|}
\prod_{\substack{1\leq k\leq n\\1\leq i\leq p}}\frac{(c_i/a_k)_{x_i}}
{(qc_i/b_k)_{x_i}}\prod_{i,k=1}^p\frac{(q^{-m_k}c_i/c_k)_{x_i}}
{(qc_i/c_k)_{x_i}}.
\end{split}\end{multline*}
\end{corollary}  

It may be interesting to note that if we choose 
$f(N)=t^N\prod_{i=1}^r(d_iq^{l_i})_N/(d_i)_N$ 
above, then  the sum in $N$ is reduced to a finite
sum by the case 
$n=1$ of Corollary \ref{c1p}. It follows that the more general series
$$\sum_{y_1,\dots,y_n=-\infty}^\infty
\frac{\De(zq^y)}{\De(z)}\prod_{i=1}^r\frac{(d_iq^{l_i})_{|y|}}{(d_i)_{|y|}}
\prod_{\substack{1\leq k\leq n\\1\leq i\leq p}}\frac{(c_iz_kq^{m_i})_{y_k}}
{(c_iz_k)_{y_k}}\prod_{i,k=1}^n\frac{(a_iz_k)_{y_k}}
{(b_iz_k)_{y_k}}\,t^{|y|} $$
may be reduced to a finite sum. The resulting identity is
awkward and we do  not write it out explicitly. However, the special case
when $p=0$ is much nicer and we state it as Corollary \ref{cn}. Note that we
have proved this
identity by first applying the case $p=0$ of Theorem \ref{t},
that is, Gustafson's identity \eqref{gi}, next Ramanujan's ${}_1\psi_1$
sum and finally the case $n=1$ of Theorem \ref{t}, which as
we have seen is due to Milne
(Corollary \ref{mc}). Thus we have actually derived  Corollary \ref{cn}  by
combining previously known results only. 

\begin{corollary}\label{cn}
For $|q^{1-|l|-n}B/A|<|t|<1$, the following identity holds:
\begin{multline*}
\sum_{y_1,\dots,y_n=-\infty}^\infty
\frac{\De(zq^y)}{\De(z)}\prod_{i=1}^r\frac{(d_iq^{l_i})_{|y|}}{(d_i)_{|y|}}
\prod_{i,k=1}^n\frac{(a_iz_k)_{y_k}}
{(b_iz_k)_{y_k}}\,t^{|y|}\\
\begin{split}&
=\frac{(AZt,q/AZt)_\infty}{(t,q^{1-n}B/At)_\infty}\frac{1}{(q^nAt/B)_{|l|}}
\prod_{i,k=1}^n\frac{(b_i/a_k,qz_k/z_i)_\infty}{(q/a_kz_i,b_iz_k)_\infty}
\prod_{i=1}^r\frac{(q^{n}d_i/BZ)_{l_i}}
{(d_i)_{l_i}}\\
&\quad\times
\sum_{x_1,\dots,x_r=0}^{l_1,\dots,l_r}\frac{\De(dq^x)}{\De(d)}\,
\left(q^{n+|l|}\frac{At}{B}\right)^{|x|}
\prod_{i=1}^r\frac{(d_i/AZ)_{x_i}}
{(q^nd_i/BZ)_{x_i}}\prod_{i,k=1}^r\frac{(q^{-l_k}d_i/d_k)_{x_i}}
{(qd_i/d_k)_{x_i}}.
\end{split}\end{multline*}
\end{corollary}

The sum on the right-hand side of Corollary \ref{cn} 
depends effectively on considerably fewer parameters than the one
on the left. We may exploit this to find a transformation formula
for the left-hand side, similarly as in Corollary \ref{st}.

\begin{corollary}\label{cmn}
Let $l_i$, be non-negative integers and $a$, $b$, $z\in\mathbb C^n$,
$\tilde a$, $\tilde b$, $\tilde z\in\mathbb C^m$, $t$, $u\in\mathbb C$ be 
parameters
with $\tilde A\tilde Z=uAZ$, $\tilde B\tilde Z=q^{m-n}uBZ$ 
and  $|q^{1-|l|-n}B/A|<|t|<1$.
Then the following transformation formula  holds:
\begin{multline*}\begin{split}&
\sum_{y_1,\dots,y_n=-\infty}^\infty
\frac{\De(zq^y)}{\De(z)}\prod_{i=1}^r\frac{(d_iq^{l_i})_{|y|}}{(d_i)_{|y|}}
\prod_{i,k=1}^n\frac{(a_iz_k)_{y_k}}
{(b_iz_k)_{y_k}}\,t^{|y|}\\
&\quad=
\frac{(AZt,q/AZt)_\infty}{(AZtu,q/AZtu)_\infty}
\prod_{i,k=1}^n\frac{(b_i/a_k,qz_k/z_i)_\infty}{(q/a_kz_i,b_iz_k)_\infty}
\prod_{i,k=1}^m\frac{(q/\tilde a_k\tilde z_i,\tilde b_i\tilde z_k)_\infty}
{(\tilde b_i/\tilde a_k,q\tilde z_k/\tilde z_i)_\infty}
\prod_{i=1}^r\frac{(ud_i)_{l_i}}{(d_i)_{l_i}}\\
&\quad\quad\times\sum_{y_1,\dots,y_m=-\infty}^\infty
\frac{\De(\tilde zq^y)}{\De(\tilde z)}\prod_{i=1}^r\frac{(ud_iq^{l_i})_{|y|}}
{(ud_i)_{|y|}}
\prod_{i,k=1}^m\frac{(\tilde a_i\tilde z_k)_{y_k}}
{(\tilde b_i\tilde z_k)_{y_k}}\,t^{|y|}.
\end{split}\end{multline*}
\end{corollary}

In the special case when $m=n$ and $\tilde a_jb_j=a_j\tilde b_j$, 
this is equivalent to Theorem 4.6 
in \cite{s}. Another interesting case is  $m=1$, when the right-hand side
reduces to a one-variable ${}_{r+1}\psi_{r+1}$, or an ${}_{r+1}\phi_r$
if we choose $\tilde b_1\tilde z_1= q$, $u=q^n/BZ$. 

\section{The case $q=1$}
\label{q1s}

In this section we consider the case of classical hypergeometric
series ($q=1$). 
Thus we will (instead of \eqref{qp}) use the standard notation
\begin{equation}\label{sn}(a)_k=\frac{\Ga(a+k)}{\Ga(a)}=
\begin{cases}a(a+1)\dotsm(a+k-1),& k\geq 0,\\
\displaystyle\frac{ 1}
{(a-1)(a-2)\dotsm(a+k)},& k<0.\end{cases}
 \end{equation}

In \cite{gu}, Gustafson proved the identities
\begin{equation}\label{5h5}\begin{split}&
\sum_{\substack{y_1,\dots,y_n=-\infty\\y_1+\dots+y_n=0}}^\infty
\frac{\De(z+y)}{\De(z)}
\prod_{i,k=1}^n\frac{(a_i+z_k)_{y_k}}
{(b_i+z_k)_{y_k}}\\
&\quad=\frac{\Ga(1+|b|-|a|-n)}{\Ga(1-|a|-|z|)\,\Ga(1+|b|+|z|-n)}
\prod_{i,k=1}^n
\frac{\Ga(1-a_k-z_i)\,\Ga(b_i+z_k)}{\Ga(b_i-a_k)\,\Ga(1+z_k-z_i)}
\end{split}\end{equation}
and 
\begin{equation}\label{2h2}\begin{split}&
\sum_{y_1,\dots,y_n=-\infty}^\infty
\frac{\De(z+y)}{\De(z)}
\prod_{\substack{1\leq k\leq n\\1\leq i\leq n+1}}
\frac{(a_i+z_k)_{y_k}}
{(b_i+z_k)_{y_k}}\\
&\quad=\frac{\Ga(|b|-|a|-n)}{\prod_{i,k=1}^{n+1}\Ga(b_i-a_k)
\prod_{i,k=1}^n\Ga(1+z_k-z_i)}
\prod_{\substack{1\leq k\leq n\\1\leq i\leq n+1}}
\Ga(1-a_i-z_k)\,\Ga(b_i+z_k),
\end{split}\end{equation}
which reduce to Dougall's ${}_5H_5$ and ${}_2H_2$ summation formulas,
when $n=2$ and $n=1$, respectively. The identity \eqref{5h5} may, 
at least formally, be obtained
 from \eqref{gi} by rescaling the parameters and letting $q$
tend to $1$.

Both these identities may be used to prove Karlsson--Minton type 
identities, exactly as in the proof of Theorem \ref{t}.
Starting with \eqref{5h5} gives the following identity, which may
also be obtained as the formal limit of Theorem \ref{t} as 
$q\rightarrow 1$.

\begin{theorem} For
$\operatorname{Re}\,(|b|-|a|)>n+|m|-1$, the following identity holds:
\begin{multline*}
\sum_{\substack{y_1,\dots,y_n=-\infty\\y_1+\dots+y_n=0}}^\infty
\frac{\De(z+y)}{\De(z)}
\prod_{\substack{1\leq k\leq n\\1\leq i\leq p}}\frac{(c_i+z_k+m_i)_{y_k}}
{(c_i+z_k)_{y_k}}\prod_{i,k=1}^n\frac{(a_i+z_k)_{y_k}}
{(b_i+z_k)_{y_k}}\\
\begin{split}&
=\frac{\Ga(1+|b|-|a|-|m|-n)}{\Ga(1-|a|-|z|-|m|)\,\Ga(1+|b|+|z|-n)}
\prod_{i,k=1}^n
\frac{\Ga(1-a_k-z_i)\,\Ga(b_i+z_k)}{\Ga(b_i-a_k)\,\Ga(1+z_k-z_i)}\\
&\quad\times
\prod_{\substack{1\leq k\leq n\\1\leq i\leq p}}\frac{(1+c_i-b_k)_{m_i}}
{(c_i+z_k)_{m_i}}
\sum_{x_1,\dots,x_p=0}^{m_1,\dots,m_p}\Bigg(\frac{\De(c+x)}{\De(c)}
\frac{(n-|b|-|z|)_{|x|}}{(1-|a|-|z|-|m|)_{|x|}}\\
&\quad\times
\prod_{\substack{1\leq k\leq n\\1\leq i\leq p}}\frac{(c_i-a_k)_{x_i}}
{(1+c_i-b_k)_{x_i}}\prod_{i,k=1}^p\frac{(c_i-c_k-m_k)_{x_i}}
{(1+c_i-c_k)_{x_i}}\Bigg).
\end{split}\end{multline*}
\end{theorem}

All the corollaries in Section \ref{s2} up to and including Corollary \ref{mnb}
have classical versions. In particular, the classical limit
of Corollary \ref{2l} is \eqref{kmg}. The classical limit of
Corollary \ref{3l} is 
\begin{equation}\label{bi}\begin{split}&
{}_{r+2}F_{r+1}\left(\begin{matrix}a,b,c_1+m_1,\dots,c_r+m_r\\
d,c_1,\dots,c_r\end{matrix}\,;1\right)\\
&\quad=(-1)^{|m|}\frac{\Ga(d)\Ga(d-a-b-|m|)}{\Ga(d-a)\Ga(d-b)}
\prod_{i=1}^r(1+c_i-d)_{m_i}\\
&\quad\quad\times\sum_{x_1,\dots,x_r=0}^{m_1,\dots,m_r}\frac{\De(c+x)}{\De(c)}
\prod_{i=1}^r\frac{(c_i-a)_{x_i}(c_i-b)_{x_i}}{(c_i)_{x_i}
(1+c_i-d)_{x_i}}
\prod_{i,k=1}^r\frac{(c_i-c_k-m_k)_{x_i}}{(1+c_i-c_k)_{x_i}}.
\end{split}\end{equation}

Using the multivariable ${}_2H_2$ sum
\eqref{2h2} one may prove the following identity.
The proof is very similar to that of 
Theorem \ref{t} but slightly easier, so we do not give
the details.

\begin{theorem} For
$\operatorname{Re}\,(|b|-|a|)>n+|m|$,  the following identity holds:
\begin{multline*}
\sum_{y_1,\dots,y_n=-\infty}^\infty
\frac{\De(z+y)}{\De(z)}
\prod_{\substack{1\leq k\leq n\\1\leq i\leq p}}\frac{(c_i+z_k+m_i)_{y_k}}
{(c_i+z_k)_{y_k}}\prod_{\substack{1\leq k\leq n\\1\leq i\leq n+1}}
\frac{(a_i+z_k)_{y_k}}{(b_i+z_k)_{y_k}}\\
\begin{split}&
=(-1)^{|m|}\frac{\Ga(|b|-|a|-|m|-n)
}{\prod_{i,k=1}^{n+1}\Ga(b_i-a_k)\prod_{i,k=1}^n\Ga(1+z_k-z_i)}
\\
&\quad\times
\prod_{\substack{1\leq k\leq n\\1\leq i\leq n+1}}
\Ga(1-a_i-z_k)\,\Ga(b_i+z_k)
\prod_{i=1}^p
\frac{\prod_{k=1}^{n+1}(1+c_i-b_k)_{m_i}}
{\prod_{k=1}^{n}(c_i+z_k)_{m_i}}\\
&\quad\times\sum_{x_1,\dots,x_p=0}^{m_1,\dots,m_p}\frac{\De(c+x)}{\De(c)}
\prod_{\substack{1\leq k\leq n+1\\1\leq i\leq p}}\frac{(c_i-a_k)_{x_i}}
{(1+c_i-b_k)_{x_i}}\prod_{i,k=1}^p\frac{(c_i-c_k-m_k)_{x_i}}
{(1+c_i-c_k)_{x_i}}.
\end{split}\end{multline*}
\end{theorem}

When $n=1$, this may be written
\begin{multline*}
\sum_{y=-\infty}^\infty\frac{(a)_y(b)_y}{(d)_y(e)_y}
\prod_{i=1}^p\frac{(c_i+m_i)_{y}}
{(c_i)_{y}}=(-1)^{|m|}\prod_{i=1}^p\frac{(1+c_i-d)_{m_i}(1+c_i-e)_{m_i}}
{(c_i)_{m_i}}\\
\begin{split}&
\times\frac{\Ga(d+e-a-b-|m|-1)\Ga(1-a)\Ga(1-b)\Ga(d)\Ga(e)}
{\Ga(d-a)\Ga(d-b)\Ga(e-a)\Ga(e-b)}
\\
&\times\sum_{x_1,\dots,x_p=0}^{m_1,\dots,m_p}\frac{\De(c+x)}{\De(c)}
\prod_{i=1}^p\frac{(c_i-a)_{x_i}(c_i-b)_{x_i}}
{(1+c_i-d)_{x_i}(1+c_i-e)_{x_i}}\prod_{i,k=1}^p\frac{(c_i-c_k-m_k)_{x_i}}
{(1+c_i-c_k)_{x_i}}.
\end{split}\end{multline*}
When $p=m_1=1$ this is a ${}_3H_3$ summation formula of Bailey
\cite{b2}, and  when $e=1$  we recover \eqref{bi}.

\end{document}